\documentclass[12pt]{article}
\usepackage{amsmath, amssymb, amsthm, amsfonts}
\usepackage{indentfirst}

\numberwithin{equation}{section}

\theoremstyle{plain}
\newtheorem{theorem}{Theorem}
\newtheorem{lemma}{Lemma}
\newtheorem{corollary}{Corollary}
\theoremstyle{definition}
\newtheorem{definition}{Definition}
\theoremstyle{remark}
\newtheorem{remark}{Remark}

\newcommand{\R}{\mathbb{R}}
\newcommand{\s}{\mathbb{S}}
\newcommand{\K}{\mathbf{k}}

\title {Local Existence for the Spatially Homogeneous Boltzmann Equation with Soft Potentials}

\author{Yong-Kum Cho \footnote{This research was supported by National Research Foundation of Korea Grant
funded by the Korean Government (\# 20120211).}}

\begin{document}

\maketitle

\begin{rmfamily}
\centerline{Department of Mathematics}

\centerline{and CAU Nonlinear PDE Center}

\centerline{College of Natural Science, Chung-Ang University}

\centerline{84 Heukseok-Ro, Dongjak-Gu, Seoul 156-756, Korea}

\medskip

\centerline{e-mail: ykcho@cau.ac.kr}

\vskip1cm

\centerline{{\bf Abstract}}
\begin{itemize}
\item[{}] We prove a local-in-time existence and uniqueness theorem for a smooth classical
solution to the spatially homogeneous Boltzmann equation with cutoff soft potentials.
Our proof is based on a series of bilinear estimates for the integrability and Sobolev regularity
of the associated collision operator. While the global-in-time existence is left inconclusive,
we give a lower bound of the maximal time of existence and a necessary condition for finite time extinction of existence.
\end{itemize}

\newpage
{\small
\begin{itemize}
\item[{}]{\bf Keywords.} Boltzmann equation, collision, contraction, cutoff, fractional integral,
Fourier transform, Hausdorff-Young type inequality, integrability, regularity, Sobolev embedding, soft potential.

\item[{}] 2010 Mathematics Subject Classification: 35Q82, 47G20, 76P05, 82B40.
\end{itemize}}

\vskip1cm

\tableofcontents

\newpage

\section{Introduction}
This paper deals with the existence and uniqueness of a smooth classical
solution to the spatially homogeneous Boltzmann equation associated with the collision kernels
classified as soft potentials.

The spatially homogeneous Boltzmann equation states
\begin{equation}\label{1.1}
\left\{\aligned &{\partial_t f(v, t) =  Q(f, f) (v, t) \quad\text{for}\quad v\in \R^d,\, t>0, }\\
& f(v, 0) = f_0(v),
\endaligned\right.
\end{equation}
which arises as a generalization of physical models in $\R^3$  for describing the behavior of a dilute gas
by its density $f$ under the simplified assumption that it depends only on the velocity $v$ and time $t$.
Here $f_0$ is a nonnegative initial datum and $Q$ stands for the bilinear integral operator
defined by
\begin{equation}\label{1.2}
Q(f, g) (v) = \int_{\mathbb{R}^d}\int_{\mathbb{S}^{d-1}} B\,
\left[f(v') g(v_*') -f(v)g(v_*)\right] d\sigma dv_*
\end{equation}
for each pair of scalar-valued functions $f, g$ on $\R^d$, where
\begin{equation}\label{1.3}
\left\{\aligned &{v' = \frac{\,v+v_*\,}{2} +
\frac{\,|v-v_*|\,}{2}\,\sigma\,,}\\
& v_*' = \frac{\,v+v_*\,}{2} - \frac{\,|v-v_*|\,}{2}\,
\sigma\,,\\
&\mathbf{k} \,\,= \frac{v-v_*}{|v-v_*|}\,,
\endaligned\right.
\end{equation}
the kernel $\,B = B(|v-v_*|, \,\mathbf{k}\cdot\sigma)\,$ is a nonnegative measurable function
and $d\sigma$ denotes the area measure on the unit sphere $\,\mathbb{S}^{d-1}\,,\,d\ge 3.\,$

In the physically relevant space $\R^3$, each pair $\,(v', v_*')\,$ represents the post-collision velocities of
two gas molecules colliding with velocities $\,(v, v_*)\,$ under binary and elastic collision dynamics. The formula (\ref{1.3})
corresponds to a parametrization of the conservation laws
\begin{equation}\label{1.4}
v+v_* = v' + v_*'\,,\quad |v|^2 + |v_*|^2 = |v'|^2 + |v_*'|^2\,.
\end{equation}
The kernel $B$ represents a physical model of collision dynamics which specifies the types of interaction potentials and impact parameters
in terms of the relative velocity $|v-v_*|$ and the deviation angle $\theta$ defined by
$\,\cos\theta = \mathbf{k}\cdot\sigma\,.$

Of primary importance is the kernel of type
\begin{equation}\label{1.5}
B = |v-v_*|^{-\lambda}\,b(\K\cdot\sigma)\qquad(-2\le\lambda<d),
\end{equation}
which generalizes the inverse-power potential model in $\R^3$.
It is classified as hard potential if $\,-2\le\lambda<0,\,$
Maxwellian if $\,\lambda=0\,$ and soft potential if $\,0<\lambda<d.$
In general, the angular part $b$ is known to be smooth
or at least bounded away from $\,\theta =0\,$ but singular at $\,\theta =0\,$ in such a way that $b$ is not integrable
over $\s^{d-1}$ as a function of $\sigma$. Upon cutting off the singularity in certain way, it is common to
assume that $\,b\in L^1(\s^{d-1})\,$, often referred to as Grad's angular cutoff assumption.

\subsection{Notions of Solution}
In the strictly classical sense, the first equation in (\ref{1.1}) means
\begin{equation}\label{1.6}
\lim_{h\to 0}\,\left|\frac{f(v, t+h) - f(v, t)}{h} - Q(f, f)(v, t)\right| = 0
\end{equation}
for each $\,(v, t)\in\R^d\times (0, \infty)\,$ with $\,Q(f,f)(v, t)<\infty.\,$

A characteristic feature of the Boltzmann equation is that nonnegativity and integrability
should be built in as a part of the definition of solution because $f$ stands for the statistical
density of a gas when $\,d=3\,.$

In consideration of these two aspects, the notion of classical solution may be defined as follows.
For simplicity, we only give a local-in-time definition on finite time interval $\,[0, T],\,T>0,\,$
without stating the detailed necessary conditions on the collision term.

\smallskip

\begin{definition}
Let $f_0$ be a nonnegative initial datum in $L^1(\R^d)$.
We say that $f$ is a \emph{classical solution} to the Boltzmann equation (\ref{1.1}) if
it is a nonnegative function satisfying (\ref{1.1}) on $\,\R^d\times [0, T],\,$
where the first equation holds in the sense of (\ref{1.6}), and $\,f(\cdot, t)\in L^1(\R^d)\,$ for all $\,t\in[0, T].\,$
\end{definition}

\smallskip

Due to the complicated nature of $Q$, it is extremely difficult to deal with
the Boltzmann equation (\ref{1.1}) in the classical sense.
In fact, we are not aware of any rigorous existence theory for the classical solutions
(we do not consider here the linearized Boltzmann
equation near equilibrium or the Boltzmann equation with small data.)
For this reason, it is common to reformulate the Boltzmann equation (\ref{1.1})
in an appropriate functional setting.

Given a Banach space $\,X\subset L^1(\R^d)\,$ with norm $\|\cdot\|_X$, the Boltzmann equation
(\ref{1.1}) in $X$ takes the form
\begin{equation}\label{1.7}
\left\{\aligned &{\partial_t f =  Q(f, f) \quad\text{in}\,\,\, X\quad\text{for}\,\,\,t\in (0, T)\,, }\\
& f(v, 0) = f_0(v)\,,
\endaligned\right.
\end{equation}
where the first equation is understood in the sense of
\begin{equation}\label{1.8}
\lim_{h\to 0}\,\left\|\frac{f(\cdot, t+h) - f(\cdot, t)}{h} - Q(f, f)(\cdot, t)\right\|_X = 0.
\end{equation}

In the case when $X$ is continuously embedded into the space of bounded continuous functions, for example,
a solution to (\ref{1.7}) becomes a classical solution,
which is a standard way of constructing classical solutions in the theory of differential equations.
Relaxing the strong differentiability of (\ref{1.8}), the integral version of (\ref{1.7}) takes the form
\begin{equation}\label{1.9}
f(v, t) = f_0(v) + \int_0^t Q(f, f)(v, s)\,ds\quad\text{for}\quad (v, t)\in\R^d\times[0, T]
\end{equation}
where the unknown $f$ is sought for in the space $\,C\left([0, T];X\right).\,$

Under certain circumstances, both notions of solutions to (\ref{1.7}) and (\ref{1.9}) could become
equivalent. Since most of existence theory for the Boltzmann equation has been studied for the integral equation (\ref{1.9})
with appropriate Banach space $X$, we only name the notion of solution for this version.

\smallskip

\begin{definition}
Given a Banach space $\,X\subset L^1(\R^d),\,$ let $f_0$ be a nonnegative initial datum in $X$.
We say that $f$ is a \emph{solution in $\,C\left([0, T];X\right)\,$} to the Boltzmann equation (\ref{1.1})
if it is a nonnegative function in  $\,C\left([0, T] ; X\right)\,$ satisfying
the integral equation (\ref{1.9}).
\end{definition}

\smallskip

There are other notions of solution to the Boltzmann equation but we shall not consider those at present.
As usual, any solution in our definitions will be said to be \emph{global} if $T$ could be arbitrary.

\subsection{Historical Backgrounds}
To the best knowledge we are aware of, the existence of a solution to the Boltzmann equation (\ref{1.1})
with kernels of type (\ref{1.5}) has been established only for the
Maxwellian and hard potentials under Grad's angular cutoff assumption on $b$.
Our overview below is based on \cite{Ar1}, \cite{CCL}, \cite{CIP}, \cite{Vi} where $L^1_\alpha(\R^d)$
denotes the space of functions $f$ on $\R^d$ with
\begin{equation}\label{1.10}
\|f\|_{L^1_\alpha} = \int_{\R^d}\,|f(v)|\,(1+|v|^2)^{\alpha/2}\,dv <+\infty\qquad(\alpha>0)
\end{equation}
and $\,T>0\,$ is arbitrary.

\begin{itemize}
\item[(1)] For the cutoff Maxwellian case,
the global existence and uniqueness for a solution in $\,C\left([0, T]; L^1(\R^d)\right)\,$
is established by D. Morgenstern (1954, \cite{Mo}) and L. Arkeryd (1972, \cite{Ar1}).
\item[(2)] For the cutoff hard potentials, the global existence of a
solution in $\,C\left([0, T]; L^1_\alpha(\R^d)\right),\,$ with $\,\alpha>2\,$ and finite entropy condition
on the initial datum, is established by L. Arkeryd (1972, \cite{Ar1}).
As it is pointed out in his papers, Arkeryd's work is based on a collection of earlier results
obtained by his predecessors including  T. Carleman (1957, \cite{Ca}), H. Grad (1965, \cite{Gr}),
A. J. Povzner (1962, \cite{Po}) who also proved the uniqueness
in the case $\,-1\le\lambda<0\,$ and E. Wild (1951, \cite{Wi}).

The uniqueness is established later by S. Mischler and B. Wennberg (1999, \cite{MW}) who also improved the solution space to
$\,C\left([0, T]; L^1_2(\R^d)\right)\,$ and removed finite entropy condition on the initial datum.
\end{itemize}

While no existence results are available in the case of soft potentials for solutions in the sense of Definition 2,
there are a set of existence results for the so-called \emph{weak solutions}. For its definition and results,
we refer to L. Arkeryd (1981, \cite{Ar2}), T. Goudon (1997, \cite{Go}), C. Villani (1998, \cite{Vi2}),
and E. Carlen, M. Carvalho and X. Lu (2009, \cite{CCL}).

A noteworthy point is that those existence results for weak solutions are obtained under non-cutoff
assumptions of the form
\begin{equation}\label{1.11}
\int_{\s^{d-1}} b(\mathbf{k}\cdot\sigma)\,( 1 - \mathbf{k}\cdot\sigma)^{j/2}\,d\sigma <+\infty
\qquad(j=1, 2).
\end{equation}
For the kernels (\ref{1.5}) of any type, where $b$ satisfies a non-cutoff condition like the above,
the existence question of a solution to the Boltzmann equation (\ref{1.1}) in the sense of Definition 2
has never been answered yet.

\subsection{Main Results}
Our primary purpose is to establish an existence theorem for a solution
to the Boltzmann equation (\ref{1.1}) with cutoff soft potentials.
In consideration of integrability and smoothness, our solution space will be
\begin{align*}
X^\alpha(\R^d) &=  L^1\cap \dot{H}^\alpha(\R^d)\\
&=\left\{\,f\in L^1(\R^d) : \| f\|_{\dot{H}^\alpha}= \left(\int_{\R^d} \,\left|\xi\right|^{2\alpha}
\bigl|\hat f(\xi)\bigr|^2\,d\xi\right)^{1/2} <\infty\,\right\}
\end{align*}
with suitable choice of $\,\alpha\ge 0\,,$ where $\hat f$ denotes the Fourier transform
\begin{equation*}
\hat f (\xi) = \int_{\R^d} e^{-2\pi i\,\xi\cdot v}\,f(v) \,dv\qquad(\xi\in\R^d)\,.
\end{equation*}

Of course, $\,\dot{H}^\alpha(\R^d)\,$ is the usual homogeneous Sobolev space of smoothing order $\alpha$
known to be a Banach space modulo polynomials. Restricting to $L^1(\R^d)$, it is clear that $X^\alpha(\R^d)$ being no longer a quotient space
is a Banach space with the maximum norm
\begin{equation}\label{1.15}
\|f\|_{X^\alpha} = \max\,\bigl(\|f\|_{L^1},\,\| f\|_{\dot{H}^\alpha}\bigr)\qquad\left(f\in X^\alpha(\R^d)\right)\,.
\end{equation}

As an alternative of Grad's angular cutoff assumption, we shall consider conditions
of type $\,b\in L^p(\s^{d-1})\,$ for some $\,1\le p\le\infty\,$ in the sense
\begin{align}\label{1.15}
\|b\|_{L^p(\mathbb{S}^{d-1})} &= \left(\int_{\mathbb{S}^{d-1}}
\left[b(\mathbf{k}\cdot\sigma)\right]^p\,d\sigma\right)^{1/p}\nonumber\\
&=\left(\left|\s^{d-2}\right|\,\int_0^\pi \left[b(\cos\theta)\right]^p\,\sin^{d-2}\theta\,d\theta
\right)^{1/p}<\infty
\end{align}
for $\,1\le p<\infty\,$ with the usual meaning for $\,p=\infty\,.$

Our principal result is the following existence and uniqueness.

\medskip

\begin{theorem} Assume that $\,B =|v-v_*|^{-\lambda}\,b(\K\cdot\sigma)\,$
where $\,b\in L^p(\mathbb{S}^{d-1})\,$ for some $\,2\le p\le\infty\,$ and $\,1/2<\lambda<d.\,$
Let $f_0$ be a nonnegative initial datum in $\,X^\alpha(\R^d) = L^1\cap \dot{H}^\alpha(\R^d)\,$ with
$\,\alpha>d/4\,$ for $\,1/2<\lambda<d/2\,$ and $\,\alpha>d/2\,$ for $\,d/2\le\lambda<d\,.$
Then there exists a finite time $T$, depending on $\,\lambda, \alpha, p, d\,$ and $\,\|b\|_{L^p(\s^{d-1})}, \,\|f_0\|_{X^\alpha}\,$ such that
the Boltzmann equation (\ref{1.1}) has a unique solution $f$ in $\,C\left([0, T] ; X^\alpha(\R^d)\right)\,$ with the following properties:
\begin{itemize}
\item[\rm{(i)}] The conservation of mass holds, that is,
\begin{equation}\label{1.16}
\int_{\R^d} f(v, t)\,dv = \int_{\R^d} f_0(v)\,dv \quad\text{for}\quad t\in [0, T].
\end{equation}
\item[\rm{(ii)}] $\,f\in C^1\left((0, T); X^\alpha(\R^d)\right)\,$ and satisfies
\begin{equation}\label{1.17}
\partial_t f =  Q(f, f) \quad\text{in}\quad X^\alpha(\R^d)\quad\text{for}\quad t\in (0, T).
\end{equation}
\end{itemize}
\end{theorem}

\medskip

In the case $\,\alpha>d/2,\,$ a basic Sobolev embedding $\,X^\alpha(\R^d)\subset C_0(\R^d)\,$
holds, where $C_0(\R^d)$ denotes the space of continuous functions $f$ on $\R^d$
which vanish at infinity, that is, $\,|f(v)|\to 0\,$ as $\,|v|\to\infty\,.$
An immediate consequence is that the solution constructed as above is indeed a classical one.
In view of its importance, we single out this result as follows.

\medskip

\begin{corollary} Assume that $\,B =|v-v_*|^{-\lambda}\,b(\K\cdot\sigma)\,$
where $\,b\in L^p(\mathbb{S}^{d-1})\,$ for some $\,2\le p\le\infty\,$ and $\,1/2<\lambda<d.\,$
Let $f_0$ be a nonnegative initial datum in $X^\alpha(\R^d)$ with $\,\alpha>d/2\,.$
Then there exists a finite time $T$, depending on $\,\lambda, \alpha, p, d\,$ and $\,\|b\|_{L^p(\s^{d-1})}, \,\|f_0\|_{X^\alpha}\,$ such that
the Boltzmann equation (\ref{1.1}) has a unique classical solution $f$ satisfying {\rm (i), (ii)} of Theorem 1
and
\begin{itemize}
\item[\rm{(iii)}] $\,f(\cdot, t)\in C_0(\R^d)\,$ for $\,t\in[0, T],\,$ $\,\,\,\partial_t f(\cdot, t)\in C_0(\R^d)\,$ for $\,t\in(0, T).\,$
\end{itemize}
\end{corollary}

\medskip

\begin{remark} Our results are local in time. While the global existence is left inconclusive, we shall discuss
in the last section about the maximal time of existence for which a lower bound will be given in terms of $\,\|f_0\|_{X^\alpha}\,.$
We also point out the following additional results:
\begin{itemize}
\item[(1)] Due to the obvious monotonicity
\begin{equation}\label{1.24}
\|b\|_{L^p(\mathbb{S}^{d-1})}\,\le\, \left|\s^{d-1}\right|^{\frac 1p - \frac 1q}\,\|b\|_{L^q(\mathbb{S}^{d-1})}\quad(1\le p<q),
\end{equation}
$L^p$ condition on $b$ gets stronger as $p$ increases. In the case $\,b\in L^p(\mathbb{S}^{d-1})\,$
with $\,1<p<2,$ there is an analogous result but valid with smaller range of $\lambda$. We postpone its
statement until the last section.
\item[(2)] In Theorem 1, it can be tracked from our proof that the solution $f$ satisfies the estimates
\begin{align}
\sup_{t\in [0, T]}\,\bigl\| f(\cdot, t)\bigr\|_{X^\alpha} &\le K\,\left\|f_0\right\|_{X^\alpha}\,,\nonumber\\
\sup_{t\in (0, T)}\,\bigl\| \partial_t f(\cdot, t)\bigr\|_{X^\alpha} &\le K\,\|b\|_{L^p(\s^{d-1})}\,\left\|f_0\right\|_{X^\alpha}^2
\end{align}
where $K$ is a constant independent of $\,f, f_0.\,$
\item[(3)] For a positive integer $\ell$, let $\,C_0^\ell(\R^d)\,$ be the space of functions $f$ on $\R^d$
such that $f$ and all of its partial derivatives up to order $\ell$  belong to $\,C_0(\R^d).\,$
In the statement of Corollary 1, if $\,\alpha >\ell + d/2\,,$ then it follows from the Sobolev embedding
$\,X^\alpha(\R^d)\subset C^\ell_0(\R^d),\,$ stated and proved in Lemma \ref{lemmaSE} below, that
the solution $f$ satisfies $\,f(\cdot, t)\in C_0^\ell(\R^d)\,$ for $\,t\in[0, T]\,$ and $\,\partial_t f(\cdot, t)\in C_0^\ell(\R^d)\,$ for $\,t\in(0, T).\,$
\end{itemize}
\end{remark}

\subsection{Methods and Outlines}
In Arkeryd's work, one of the essential ingredients is the bilinear estimate
\begin{equation}\label{1.19}
\left\| Q(f, g)\right\|_{L^1}\,\le\, 2\|b\|_{L^1(\s^{d-1})}\,\|f\|_{L^1}\,\|g\|_{L^1}
\end{equation}
for the cutoff Maxwellian case. Indeed, it implies that $Q$ is locally Lipschitz so that the local existence
and uniqueness, except nonnegativity, follows right away by the Banach contraction mapping principle.
Due to the conservation of mass, the global existence follows by repeating
the same arguments on the next time interval of equal length iteratively. For the cutoff hard-potential case, the method of proof
relies on a fine scheme of extracting a weakly convergent subsequence in $L^1_\alpha$ from the sequence of solutions
obtained by truncating the potential part and applying the already-established existence theorem for the Maxwellian case.

In the same spirit as in Arkeryd's work for the Maxwellian case, our work will be based on
a bilinear estimate of $Q$ in the form
\begin{equation}\label{1.20}
\left\| Q(f, g)\right\|_{X^\alpha}\,\le\, K\|b\|_{L^p(\s^{d-1})}\,\|f\|_{X^\alpha}\,\|g\|_{X^\alpha}\,,
\end{equation}
which yields the local existence and uniqueness at once except nonnegativity.
For nonnegativity, we shall adopt the method of X. Lu and Y. Zhang (\cite{LZ}) for which the key idea is
to prove that the negative part $\,(-f)^+\,$ of the solution $f$ obtained has zero mass
\begin{equation}\label{1.21}
\int_{\R^d}\left( -f(v, t)\right)^+\,dv =0 \quad\text{for}\quad t\in [0, T].
\end{equation}

As for the bilinear estimate (\ref{1.20}), we shall investigate
the gain part $Q^+$ and the loss part $Q^-$
separately, where $\,Q= Q^+ - Q^-\,$ and
\begin{align}\label{1.21}
Q^+(f, g) (v) &= \int_{\mathbb{R}^d}\int_{\mathbb{S}^{d-1}} B\, f(v')\, g(v_*')\,d\sigma dv_*\,,\\
Q^-(f, g) (v) &= f(v)\int_{\mathbb{R}^d}\int_{\mathbb{S}^{d-1}} B\, g(v_*)\,d\sigma dv_*\,.
\end{align}
While the integrability is a simple consequence of the pointwise
boundedness of fractional integrals on $X^\alpha(\R^d)$, the regularity
is subtle. In fact, it turns out that $Q^+$ has regularity-gaining properties,
whereas $Q^-$ has regularity-preserving properties.
By using the monotone embedding of the $X^\alpha(\R^d)$,
we shall obtain regularity estimates for $Q$ in a collective manner.

Besides the global existence, there are two other deficiencies in our work. One is that we were not able to
cover the range $\,0<\lambda\le 1/2\,.$ Another is that
our $L^p$ cutoff assumptions on $b$ are restricted to $\,p>1\,.$
As we shall see below, those deficiencies are partly due to the fact that our regularity analysis on $Q^+$ relies on
exploiting certain decay estimates of the Fourier multiplier arising in its Fourier transform which are case sensitive
to the potential order $\lambda$ and the type of cutoff condition on $b$.

Throughout this paper, most of multiplicative constants appearing
in our estimates will be computable but we shall not care for the best possible ones. For the sake of
reminding or simplifying, we shall use the notation $L^1\cap\dot{H}^\alpha(\R^d)$ and $X^\alpha(\R^d)$ interchangeably.

\section{Function Spaces and Fractional Integrals}
As it is demonstrated in our recent work \cite{Cho},
the $\,X^\alpha(\R^d) =L^1\cap\dot{H}^\alpha(\R^d)\,$ provide an ideal framework for studying the integrability
and regularity of the collision operator $Q$ associated with soft potentials.

One of the main advantages is that Hausdorff-Young type inequalities and Sobolev type
inequalities are available in the following logarithmic convexity forms
(see Lemma 2.1 and Lemma 2.2, \cite{Cho}).

\medskip

\begin{lemma}\label{lemmaHY1} {\rm (Hausdorff-Young type inequality)}
For $\,1\le p\le\infty\,,$ let $\alpha$ be any number satisfying $\,\alpha>(1/p - 1/2)d\,$ when $\,1\le p<2\,$ and $\,\alpha\ge 0\,$ when $\,2\le p\le\infty\,.$ If
$\,f\in L^1\cap\dot{H}^\alpha(\R^d)\,,$ then $\,\hat f\in L^p(\R^d)\,$ and
\begin{equation}\label{HY2}
\bigl\|\hat f\bigr\|_{L^p}\,\le\,C\,\|f\|_{L^1}^{1-\mu}\,\|f\|_{\dot{H}^\alpha}^\mu\quad\text{with}\quad \mu = \frac{2d}{p(d+2\alpha)}\in [0,1],
\end{equation}
where $C$ is a constant depending only on $\,\alpha, d, p\,.$
\end{lemma}

\medskip

\begin{lemma}\label{lemmaHY2} {\rm (Sobolev type inequality)}
For $\,1\le p\le\infty\,,$ let $\alpha$ be any number satisfying $\,\alpha>(1/2-1/p)d\,$ when $\,2<p\le\infty\,$ and $\,\alpha\ge 0\,$ when $\,1\le p\le 2\,.$ If
$\,f\in L^1\cap\dot{H}^\alpha(\R^d)\,,$ then $\,f\in L^p(\R^d)\,$ and
\begin{equation}\label{HY3}
\bigl\|f\bigr\|_{L^p}\,\le\,C\,\|f\|_{L^1}^{1-\nu}\,\|f\|_{\dot{H}^\alpha}^\nu\quad\text{with}\quad \nu = \frac{2d}{d+2\alpha}
\left(1-\frac 1p\right)\in [0,1],
\end{equation}
where $C$ is a constant depending only on $\,\alpha, d, p\,.$
\end{lemma}

An application of Lemma \ref{lemmaHY1} with the particular instance $\,p=1\,$ yields
the following Sobolev embedding properties which
will be used in lifting any solution in $X^\alpha(\R^d)$ to a classical solution.

\medskip

\begin{lemma}\label{lemmaSE}
Let $\,f\in X^\alpha(\R^d).$
\begin{itemize}
\item[{\rm(i)}] If $\,\alpha>d/2,\,$ then $\,f\in C_0(\R^d)\,$ with
\begin{equation}
\sup_{v\in\R^d}\,|f(v)|\,\le\,K\,\|f\|_{X^\alpha}
\end{equation}
for some constant $K$ independent of $f$.
\item[\rm{(ii)}] If $\,\alpha>\ell + d/2\,$ with a positive integer $\ell$, then $\,f\in C_0^\ell(\R^d)\,$ with
\begin{equation}
\sum_{|m|\le\ell}\left[\sup_{v\in\R^d}\,\bigl|\partial^m f(v)\bigr|\right]\,\le\,K\,\|f\|_{X^\alpha}
\end{equation}
for some constant $K$ independent of $f$, where $\,m=(m_1, \cdots, m_d)\in\mathbb{Z}_+^d\,$ denotes a multi-index and
$\,|m| = m_1 + \cdots + m_d\,.$
\end{itemize}
\end{lemma}

\smallskip

\begin{proof}
For (i), Lemma \ref{lemmaHY1} shows $\,\hat f\in L^1(\R^d)\,$ and hence
\begin{equation}\label{SE1}
f(v) = \int_{\R^d} e^{2\pi i\,v\cdot\xi}\,\hat f(\xi) \,d\xi\qquad(v\in\R^d)
\end{equation}
by the Fourier inversion theorem. In view of the Riemann-Lebesgue lemma, we conclude $\,f\in C_0(\R^d)\,$
and the stated estimate is a simple consequence of (\ref{HY2}) and AM-GM inequality.

For part (ii), observe that
\begin{equation}\label{SE2}
\int_{\R^d}|\xi|^{|m|}\,\bigl|\hat f(\xi)\bigr|\,d\xi \,\le\,C\,\|f\|_{L^1}^{1-\delta}\,\|f\|_{\dot{H}^\alpha}^{\delta}\,,
\quad \delta = \frac{2(d +|m|)}{d+2\alpha}\,,
\end{equation}
which can be verified easily by the method of splitting and optimizing
(see the proof of Lemma \ref{lemmaHY3} below). Differentiating (\ref{SE1}) under the integral sign,
the result follows by applying (i) repeatedly with (\ref{SE2}).
\end{proof}

\medskip
The following monotonicity will be used in combining regularity estimates of $Q^+$
for which the regularity orders are different.

\medskip

\begin{lemma}\label{lemmaHY3}
Let $\,0\le\beta\le\alpha\,.$
If $\,f\in L^1\cap\dot{H}^\alpha(\R^d),\,$ then $\,f\in \dot{H}^\beta(\R^d)\,$ and there exists a constant $C$ depending only on $\,\alpha, \beta, d\,$ such that
\begin{equation}\label{HY4}
\|f\|_{\dot{H}^\beta}\,\le\,C\, \|f\|_{L^1}^{1-\kappa}\,\|f\|_{\dot{H}^\alpha}^\kappa
\quad\text{with}\quad \kappa = \frac{d + 2\beta}{d+2\alpha}\in (0, 1].
\end{equation}
As a consequence, the embedding $\,X^\alpha(\R^d)\subset X^\beta(\R^d)\,$ holds with
\begin{equation}
\|f\|_{X^\beta}\,\le\,C\, \|f\|_{X^\alpha}\quad\text{for}\quad f\in X^\alpha(\R^d)\,,
\end{equation}
where the constant $C$ is independent of $f$.
\end{lemma}

\smallskip

\begin{proof} For any $\,\rho>0\,,$ if we split the integral into two parts by writing $\,\R^d = \{|\xi|\le\rho\}\cup \{|\xi| >\rho\}\,,$
then it is trivial to observe
$$
\int_{\R^d}|\xi|^{2\beta}\,\bigl|\hat{f}(\xi)\bigr|^2\,d\xi\,\le\,C\left\{\rho^{2\beta+d}\,\|f\|_{L^1}^2 + \rho^{2\beta -2\alpha}\,\|f\|_{\dot{H}^\beta}^2\right\}
$$
with the constant $C$ depending only on $\,\beta, d\,.$ Optimizing in $\rho$ and taking the square root, we obtain the desired conclusion (\ref{HY4}).
\end{proof}

\medskip

Another important advantage of the $X^\alpha(\R^d)$ is that fractional integrals are uniformly bounded in the pointwise sense once $\alpha$
were chosen suitably, which will be used in establishing the integrability of $Q$.
For later use, we label the multiplicative constant in the corresponding estimate.

\medskip

\begin{lemma}\label{lemmaF1} For $\,0<\lambda<d\,,$ put
\begin{equation}\label{F1}
J_\lambda(f) = \sup_{v\in\R^d}\,\int_{\R^d}\,|v- v_*|^{-\lambda}\,| f(v_*)|\,d v_* .
\end{equation}
Let $\alpha$ be any number satisfying $\,\alpha\ge 0\,$ when $\,0<\lambda<d/2\,$ and $\,\alpha>\lambda - d/2\,$ when $\,d/2\le\lambda<d.\,$
If $\,f\in L^1\cap\dot{H}^\alpha(\R^d)\,,$ then there exists a constant $C_J$ depending only on $\,\alpha, \lambda, d\,$ such that
\begin{equation}\label{F2}
J_\lambda(f)\,\le\,C_J\,\|f\|_{L^1}^{1-\theta}\,\|f\|_{\dot{H}^\alpha}^\theta\quad \text{with}\quad \theta = \frac{2\lambda}{d+2\alpha}\,\in\,(0,1).
\end{equation}
\end{lemma}

\smallskip

\begin{proof} Fix $\,v\in\R^d.\,$ For any $\,\rho>0\,$ and $\,1<p<d/\lambda,\,$
we split the integral into two parts by writing $\,\R^d = \{|v_* -v|\le\rho\}\cup \{|v_* -v| >\rho\}\,$
and apply H\"older's inequality to estimate
\begin{align*}
\int_{\R^d}\,|v- v_*|^{-\lambda}\, |f(v_*)|\, d v_* \,
\le\, C\,\left\{ \rho^{-\lambda +d/p}\,\|f\|_{L^{p'}} + \rho^{-\lambda}\,\|f\|_{L^1}\right\}
\end{align*}
where $C$ depends only on $\,\lambda, p, d\,$ and $\,1/p \,+ \,1/p' = 1\,.$
Optimizing in $\rho$ gives
\begin{equation*}
\int_{\R^d}\,|v- v_*|^{-\lambda}\, |f(v_*)|\, d v_* \, \le\, C\, \|f\|_{L^1}^{1- p\lambda/d}\,\| f\|_{L^{p'}}^{p\lambda/d}\,.
\end{equation*}

An application of Sobolev type inequality (\ref{HY3}) yields
\begin{equation*}
\|f\|_{L^{p'}}\, \le\, C\, \|f\|_{L^1}^{1- \delta}\,\|f\|_{\dot{H}^\alpha}^{\delta}\quad\text{with}\quad  \delta = \frac{2d}{(d+2\alpha)p}
\end{equation*}
provided $\,\alpha\ge 0\,$ when $\,p\ge 2\,$ and $\,\alpha>(1/p -1/2)d\,$ when $\,1<p<2,\,$ where
$C$ depends on $\,\alpha, \lambda, p, d\,$ in this occasion. If we choose $\,p= 2\,$ in the case $\,0<\lambda<d/2\,$
and $\,p= d/2\lambda\,$ in the case $\,\lambda\ge d/2\,,$ for instance, then the above condition is trivially verified.
In view of the relation $\,\delta\,p\lambda/d = \theta,\,$
inserting the last estimate and simplifying the exponents with the above choice of $p$, we conclude
\begin{equation*}
\int_{\R^d}\,|v- v_*|^{-\lambda}\, |f(v_*)|\, d v_* \, \le\, C\,\|f\|_{L^1}^{1-\theta}\,\|f\|_{\dot{H}^\alpha}^\theta,
\end{equation*}
where $C$ depends only on $\,\alpha, \lambda, d,\,$ and the proof is complete.
\end{proof}

\medskip

\begin{remark} A counterpart of the estimate (\ref{F2}) is given by
\begin{equation}\label{F-3}
\sup_{\xi\in\R^d}\,\int_{\R^d}\,|\zeta|^{-d+\lambda}\,\bigl|\hat f(\xi -\zeta)\bigr|\,d\zeta
\,\le\,C\,\|f\|_{L^1}^{1-\theta}\,\|f\|_{\dot{H}^\alpha}^\theta
\end{equation}
with the same $\theta$, which is proved in Lemma 3.1 of \cite{Cho}. A well-known Fourier
representation of the fractional integral states
\begin{align}\label{F3}
&\int_{\R^d}\,|v- v_*|^{-\lambda}\, f(v_*)\, d v_* = \lambda_d\int_{\R^d} e^{2\pi i\, v\cdot\zeta}\,
|\zeta|^{-d+\lambda}\,\hat f(\zeta)\,d\zeta,\nonumber\\
&\text{where}\quad \lambda_d =  \pi^{\frac{-d+\lambda}{2}}\,\Gamma\left(\frac{d-\lambda}{2}\right)\Big/\Gamma\left(\frac{\lambda}{2}\right),
\end{align}
for each $\,v\in\R^d\,$ and $\,0<\lambda<d\,$ (see \cite{St} for example). The estimate (\ref{F2}) for $J_\lambda(f)$ is a simple consequence of (\ref{F-3}) in the case
$\,f\ge 0.$ Throughout this paper, $\lambda_d$ will stand for the constant defined as in (\ref{F3}).
\end{remark}

\section{Integrability and Regularity}
Our ultimate goal in this section is to establish bilinear estimates for the collision
operator $Q$ with cutoff soft potentials in the aforementioned form (\ref{1.20}). To this end,
we shall deal with the integrability and the regularity of $Q$ separately. As it turns out to be
quite simple to obtain the integrability estimate, we shall focus on studying the regularity of $Q$ for which a series of
bilinear estimates for the gain part $Q^+$ and the loss part $Q^-$
will be proved case by case and combined to yield the desired estimate.

\subsection{Integrability of $Q$}
Since the collision operators with soft potentials are closely related with fractional integrals,
it is a simple matter to establish its integrability once the pointwise
estimates for fractional integrals were available.

\medskip

\begin{theorem}\label{theoremI1} Assume that $\,B = |v-v_*|^{-\lambda}\,b\left(\mathbf{k}\cdot\sigma\right)\,$
where $\,b\in L^1(\mathbb{S}^{d-1})\,$ and $\,0<\lambda<d\,$.
If $\,f, g\in L^1\cap\dot{H}^\alpha(\R^d)\,$ with $\,\alpha\ge 0\,$ for $\,0<\lambda<d/2\,$ and $\,\alpha>\lambda - d/2\,$ for $\,d/2\le\lambda<d,\,$
then
\begin{align}\label{I1}
&\left\|Q(f,g)\right\|_{L^1}\,\le\, 2C_J\,\|b\|_{L^1(\mathbb{S}^{d-1})}\,\|f\|_{L^1}^{1-\theta}\,\|f\|_{\dot{H}^\alpha}^\theta\,\|g\|_{L^1}\\
&\text{with}\quad \theta=\frac{2\lambda}{d+2\alpha}\in(0, 1),\nonumber
\end{align}
where $C_J$ denotes the same constant as in (\ref{F2}) of Lemma \ref{lemmaF1}.
\end{theorem}

\smallskip

\begin{proof}
By using the fact that the transformation $\,(v, v_*, \sigma)\to (v', v_*', \mathbf{k})\,$ has unit Jacobian and applying
Lemma \ref{lemmaF1}, we deduce
\begin{align}\label{I2}
\int\left|Q(f,g)(v)\right|dv &\le \int\left|Q^+(f,g)(v)\right|dv  + \int\left|Q^-(f,g)(v)\right|dv \nonumber\\
&\le 2\,\|b\|_{L^1(\mathbb{S}^{d-1})}\,\iint |v-v_*|^{-\lambda}|f(v)||g(v^*)|\,dv_* dv\nonumber\\
&\le 2\,\|b\|_{L^1(\mathbb{S}^{d-1})}\,J_\lambda(f)\,\|g\|_{L^1}
\end{align}
and we obtain the result by invoking Lemma \ref{lemmaF1} for an estimate of $J_\lambda(f)$.
\end{proof}

\medskip
What we aim to establish for the integrability of $Q$ is the following:

\medskip

\begin{corollary}\label{corollaryI1} Assume that $\,B = |v-v_*|^{-\lambda}\,b\left(\mathbf{k}\cdot\sigma\right)\,$
where $\,b\in L^p(\mathbb{S}^{d-1})\,$ for some $\,1\le p\le\infty\,$ and $\,0<\lambda<d\,$.
If $\,f, g\in X^\alpha(\R^d)\,$ with
$\,\alpha\ge 0\,$ for $\,0<\lambda<d/2\,$ and $\,\alpha>\lambda - d/2\,$ for $\,d/2\le\lambda<d,\,$
then there exists a constant $C_I$ depending only on $\,\lambda, \alpha, p, d\,$ such that
\begin{equation}\label{I2}
\left\|Q(f,g)\right\|_{L^1}\,\le\, C_I\,\|b\|_{L^p(\mathbb{S}^{d-1})}\,\|f\|_{X^\alpha}\,\|g\|_{X^\alpha}\,.
\end{equation}
\end{corollary}

\smallskip

\begin{proof} Apply Theorem \ref{theoremI1} and put $\,C_I = 2 C_J\,|\s^{d-1}|^{1-1/p}\,.$
\end{proof}

\subsection{Regularity of $Q^+$}
By the preceding result of integrability, it is possible to take the Fourier transform of the gain operator $Q^+$
defined by
$$Q^+(f, g)(v) = \int_{\R^d}\int_{\s^{d-1}} |v-v_*|^{-\lambda}\,b\left(\mathbf{k}\cdot\sigma\right)\,f(v')\, g(v_*')\, d\sigma dv_*$$
in the pointwise sense for each fixed functions $\,f, g\in X^\alpha(\R^d)\,$ if $\alpha$ satisfies the size condition
described as in Corollary \ref{corollaryI1}.

In fact, if we make use of the Fourier representations (\ref{F3}) of fractional integrals and the idea of Bobylev (\cite{Bo})
for evaluating surface integrals by exchanging unit vectors, it is straightforward to derive
\begin{align}\label{N0}
[Q^+(f,g)]\,\widehat{\,}\,(2\xi) &= \lambda_d\int_{\R^d} m(\xi, \zeta)
\,\hat f(\xi + \zeta) \hat g(\xi - \zeta)\, d\zeta\,,\nonumber\\
m(\xi, \zeta) &= \int_{\mathbb{S}^{d-1}}b\biggl(\frac{\xi}{|\xi|} \cdot\sigma\biggr)
\bigl|\zeta -|\xi|\sigma\bigr|^{-d +\lambda} \,d\sigma\,.
\end{align}

In our recent work \cite{Cho}, this Fourier transform formula was used in extending the celebrated theorem of P.-L. Lions (\cite{Li})
on the $\,(d-1)/2\,$-regularity-gaining property of $Q^+$ to the case of cutoff soft potentials.
As it is still relevant to the present work, we quote one of those regularity results which
corresponds to the particular instance $\,\beta =0\,$ of Theorem 4.5, \cite{Cho}.

\medskip

\begin{theorem}\label{theoremR1} Assume that $\,B = |v-v_*|^{-\lambda}\,b\left(\mathbf{k}\cdot\sigma\right)\,$
where $\,b\in L^p(\mathbb{S}^{d-1})\,$ for some $\,1<p\le\infty\,$ and $\,1+(d-1)/p<\lambda<d\,$.
If $\,f, g\,\in\,L^1\cap\dot{H}^\alpha(\R^d)\,$ with any nonnegative $\alpha$ satisfying
$\,\alpha>\lambda -d/2\,,$ then there exists a constant $C$ depending only on $\,\lambda, \alpha, p, d\,$
such that
\begin{align}\label{R1}
\left\|Q^+(f,g)\right\|_{\dot{H}^{\alpha}}&\le C\|b\|_{L^p(\mathbb{S}^{d-1})}\biggl\{\|f\|_{\dot{H}^\alpha} \|g\|_{L^1}^{1-\theta}\|g\|_{\dot{H}^\alpha}^\theta
+\|g\|_{\dot{H}^\alpha}\|f\|_{L^1}^{1-\theta}\|f\|_{\dot{H}^\alpha}^\theta \biggr\}\nonumber\\
&\text{where}\quad \theta = \frac{2\lambda}{d+2\alpha}\in (0, 1)\,.
\end{align}
\end{theorem}

\medskip

The proof is based on obtaining the decay rates of $m(\xi, \zeta)$
with the aid of the following elementary estimate for surface integrals (see Lemma 2.4, \cite{Cho})
and dominating the Fourier transform of $Q^+$ by certain bilinear fractional integral
whose $L^p$ behaviors are known (see \cite{KS}).

\medskip

\begin{lemma}\label{lemmaS} For $\,\xi, \zeta\in\R^d\,$ with $\,d\ge 3\,,$ if $\,a<d-1\,,$ then
\begin{equation}\label{S1}
\int_{\mathbb{S}^{d-1}}\bigl|\zeta - |\xi|\sigma\bigr|^{-a}\,d\sigma\, \le\,C\,
\bigl(|\xi| + |\zeta|\bigr)^{-a}
\end{equation}
for some constant $C$ depending only on $\,a, d\,.$
\end{lemma}

\medskip

To deal with other values of $\lambda$, we prove the following:

\medskip

\begin{theorem}\label{theoremR2} Assume that $\,B = |v-v_*|^{-\lambda}\,b\left(\mathbf{k}\cdot\sigma\right)\,$
where $\,b\in L^p(\mathbb{S}^{d-1})\,$ for some $\,1<p\le 2\,$ and
\begin{equation}\label{N1}
\frac dp - (d-1)\left(1-\frac 1p\right)<\lambda\le\frac dp\,.
\end{equation}
If $\,f, g\,\in\,L^1\cap\dot{H}^\alpha(\R^d)\,$ with $\alpha$ satisfying
\begin{equation}\label{N-1}
\alpha>\max\left\{ \frac 1p\left(1-\frac 1p\right)d, \, \left(\frac 1p -\frac 12\right) d\right\}\,,
\end{equation}
then there exists a constant $C$ depending only on $\,\alpha, \lambda, p, d\,$
such that
\begin{align}\label{N2}
&\left\|Q^+(f,g)\right\|_{\dot{H}^{\alpha + \frac dp -\lambda}}\nonumber\\
&\qquad\le\, C\,\|b\|_{L^p(\mathbb{S}^{d-1})}\biggl\{\|f\|_{\dot{H}^\alpha} \|g\|_{L^1}^{1-\mu}\|g\|_{\dot{H}^\alpha}^\mu
+\|g\|_{\dot{H}^\alpha}\|f\|_{L^1}^{1-\mu}\|f\|_{\dot{H}^\alpha}^\mu \biggr\}\nonumber\\
&\,\text{where}\quad \mu = \frac{2d}{p(d+2\alpha)}\in (0,1)\,.
\end{align}
\end{theorem}

\smallskip

\begin{proof} For simplicity, put
$\,\widehat{Q^+}\,(\xi) = \left[Q^+(f, g)\right]\widehat{\,}\,\,(2\xi)\,.$
Let $\,1/q=1-1/p\,$ and fix momentarily a number $\beta$ satisfying
\begin{equation}\label{N3}
\max\left\{\frac dp, \,\lambda\right\} <\beta < \min\left\{\lambda + \frac{d-1}{q}, \,d\right\}.
\end{equation}
By using the beta integral identity and Lemma \ref{lemmaS}, we deduce
\begin{align}
\left(\int_{\R^d} \bigl|\zeta -|\xi|\sigma\bigr|^{-q(d -\beta)}|\zeta|^{-\beta}\,d\zeta\right)^{1/q}\,&\le\,C\,
|\xi|^{-\left(\frac{d-\beta}{p}\right)}\,,\label{N4}\\
\left(\int_{\s^{d-1}} \bigl|\zeta -|\xi|\sigma\bigr|^{q(\lambda-\beta)}\,d\sigma\right)^{1/q} \,&\le\,C\,(|\xi| +|\zeta|)^{\lambda -\beta}\,.\label{N5}
\end{align}

The estimate (\ref{N4}) and H\"older's inequality yield
\begin{align}\label{N6}
&\quad\int \bigl|\zeta -|\xi|\sigma\bigr|^{-d +\lambda} \bigl|\hat f(\xi +\zeta) \hat g(\xi -\zeta)\bigr|\,d\zeta\,\le\,C\,|\xi|^{-\left(\frac{d-\beta}{p}\right)}
\nonumber\\
&\qquad\quad\times\quad
\left(\int\bigl|\zeta -|\xi|\sigma\bigr|^{p(\lambda-\beta)}|\zeta|^{\frac{p\beta}{q}} \bigl|\hat f(\xi +\zeta) \hat g(\xi -\zeta)\bigr|^p\,d\zeta\right)^{1/p}.
\end{align}
With the aid of the estimates (\ref{N5}), (\ref{N6}),  if we apply H\"older's inequality first and Minkowski's integral inequality subsequently,
valid due to $\,q/p \ge 1\,,$ then it is straightforward to derive
\begin{align}\label{N7}
&\left|\widehat{Q^+}\,(\xi)\right| \le C\, \|b\|_{L^p(\s^{d-1})}\,|\xi|^{-\left(\frac{d-\beta}{p}\right)}\quad\times\nonumber\\
&\quad\quad \left(\int (|\xi| +|\zeta|)^{p(\lambda-\beta)}|\zeta|^{\frac{p\beta}{q}} \bigl|\hat f(\xi +\zeta) \hat g(\xi -\zeta)\bigr|^p\,d\zeta\right)^{1/p}.
\end{align}

For any $\,\delta\ge 0\,,$ we observe that
\begin{align*}
(|\xi| +|\zeta|)^{\lambda-\beta}\,|\zeta|^{\frac{\beta}{q}} \le C\,|\xi|^{\lambda-\beta -\delta}\,\left(|\xi +\zeta|^{\delta +\frac{\beta}{q}} +
|\xi -\zeta|^{\delta + \frac{\beta}{q}}\right)
\end{align*}
with a constant $C$ depending only on the exponents. By using this elementary inequality, we deduce from (\ref{N7})
\begin{align*}
&\qquad\left|\widehat{Q^+}\,(\xi)\right| \le C\, \|b\|_{L^p(\s^{d-1})}\,|\xi|^{-\left(\frac{d-\beta}{p}\right) +\lambda -\beta -\delta}\quad\times\nonumber\\
&\biggl\{ \left(\int |\widehat F(\xi +\zeta)\hat g(\xi -\zeta)\bigr|^p d\zeta\right)^{1/p} +
\left(\int |\hat f(\xi +\zeta)\widehat G(\xi -\zeta)\bigr|^p d\zeta\right)^{1/p}\biggr\}
\end{align*}
where we put $\,\widehat F (\xi ) = |\xi|^{\delta + \frac{\beta}{q}} \,\hat f(\xi),\, \widehat G(\xi ) = |\xi|^{\delta + \frac{\beta}{q}}\, \hat g(\xi)\,.$
It follows from the ordinary Minkowski inequality that
\begin{align}
&\left\||\xi|^{\left(\frac{d-\beta}{p}\right) -\lambda +\beta +\delta}\,\widehat{Q^+}\,(\xi)\right\|_{L^2(d\xi)}
\le\, C\,\|b\|_{L^p(\s^{d-1})}\,\bigl\{{\rm(I)} + {\rm (II)}\bigr\}\label{N8}\\
\text{where}&\qquad\quad {\rm(I)}= \biggl\|\left(\int |\widehat F(\xi +\zeta) \,\hat g(\xi -\zeta)\bigr|^p\,d\zeta\right)^{1/p}\biggr\|_{L^2(d\xi)}\,,\nonumber\\
&\qquad\quad {\rm(II)}= \biggl\|\left(\int |\hat f(\xi +\zeta) \,\widehat G(\xi -\zeta)\bigr|^p\,d\zeta\right)^{1/p}\biggr\|_{L^2(d\xi)}\,.\nonumber
\end{align}

Since $\,2/p\ge 1\,,$ we may apply Minkowski's integral inequality to estimate
\begin{align*}
{\rm(I)} &= \biggl\|\left(\int |\widehat F(2\xi -\zeta)\, \hat g(\zeta)\bigr|^p\,d\zeta\right)^{1/p}\biggr\|_{L^2(d\xi)}\\
&\le \biggl[\int\left(\int |\widehat F(2\xi -\zeta)\, \hat g(\zeta)\bigr|^2\,d\xi\right)^{p/2} d\zeta\biggr]^{1/p}\\
& = 2^{-d/2}\,\bigl\|\widehat F\bigr\|_{L^2}\,\left\|\hat g\right\|_{L^p}\,.
\end{align*}
Interchanging the roles of $\,f, g\,$ and $\,F, G\,$, we also get
$$ {\rm(II)}\,\le\, 2^{-d/2}\,\bigl\|\widehat G\bigr\|_{L^2}\,\left\|\hat f\right\|_{L^p}\,.$$
Upon setting $\,\alpha = \delta + \beta/q\,$ and simplifying, we conclude from (\ref{N8})
\begin{equation}\label{N9}
\left\|Q^+\right\|_{\dot{H}^{\alpha + \frac dp -\lambda}}\,\le\,
C\,\|b\|_{L^p(\s^{d-1})}\,\left\{\|f\|_{\dot H^\alpha}\,\left\|\hat g\right\|_{L^p} + \|g\|_{\dot H^\alpha}\,\bigl\|\hat f\bigr\|_{L^p}\right\}\,,
\end{equation}
which yields immediately the desired estimate (\ref{N2}) in view of the Hausdorff-Young type inequality (\ref{HY2}), provided
$\,\alpha>(1/p -1/2) d\,.$

It remains to check the necessary conditions (\ref{N1}), (\ref{N-1}). Concerning $\lambda$,
the choice condition (\ref{N3}) on $\beta$ requires
$$
\frac dp<\lambda + \frac{d-1}{q}\quad\text{or equivalently}\quad \frac dp - (d-1)\left(1-\frac 1p\right)<\lambda
$$
for such a $\beta$ would not exist otherwise. The requirement $\,\lambda\le d/p\,$ is needed to ensure the
regularity gaining property of $Q^+$. On the other hand, since $\,\alpha>\beta/q\,$, the choice condition (\ref{N3}) on $\beta$ implies
$$\alpha>\max\left\{\left(1-\frac 1p\right)\frac dp, \,\left(1-\frac 1p\right)\lambda\right\}
= \left(1-\frac 1p\right)\frac dp\,.$$
Comparing with $\,\alpha > (1/p -1/2) d\,,$ it is easy to see that (\ref{N-1}) comes up necessarily
for the validity of our arguments. The proof is now complete.
\end{proof}

\medskip

\begin{remark} Under the assumption $\,b\in L^p(\s^{d-1})\,$ with
$\,2<p\le\infty,\,$ if we interchange the roles played by $p, q$ in the above proof and proceed in the
same way, then it is a routine matter to find that $Q^+$ gains regularity of order
$\,(1-1/p)d -\lambda\,$ for $\lambda$ in the range
$$1-\frac 1p<\lambda\le \left(1-\frac 1p \right) d$$
if $\alpha$ satisfies certain size condition similar to (\ref{N-1}).
As we are only concerned about the maximum possible range of $\lambda$
for which regularity gaining or preserving property of $Q^+$ continues to hold, we shall not need
such a result.
\end{remark}

\medskip

Suppose now that $\,b\in L^2(\s^{d-1}).\,$ By Theorem \ref{theoremR1} and Theorem \ref{theoremR2},
$Q^+$ has a regularity-gaining or preserving property on $X^\alpha(\R^d)$ for
$$\lambda\in \left(\frac 12, \,\frac d2\right]\,\cup\,\left(\frac{d+1}{2}, \,d \right)\,$$
with $\,\alpha>d/4\,$ or $\,\alpha>\lambda - d/2\,$ on each interval, respectively.
If we recall $\,L^2(\s^{d-1})\subset L^p(\s^{d-1})\,$ for all $\,1<p\le 2\,$
and apply Theorem \ref{theoremR2} repeatedly, we may extend such a result to one of the remaining intervals.

To see this, let us write the $\lambda$-interval of validity in Theorem \ref{theoremR2}
as
$$\Phi(p) = \frac dp - (d-1)\left(1 -\frac 1p\right)\,,\,\, I(p) = \left(\Phi(p),\,\frac dp\right].$$
As $p$ increases from $1$ to $2$, the end-point function $\Phi$ decreases from $d$ to $\frac 12$ and the interval length of $I(p)$
increases from $0$ to $\frac{d-1}{2}$. Thus
\begin{equation*}\label{NN3}
\bigcup_{1<p\le 2}\,I(p) = \left(\frac 12,\,d\right),\,
\end{equation*}
which is the maximum possible range of $\lambda$ for which $Q^+$ has a regularity gaining property
obtainable from Theorem \ref{theoremR2}.

A simple inspection shows, for instance,
$$\,\left(\frac d2, \,\frac {d+1}{2}\right]\subset I(p)\quad\text{with}\quad
p=\frac{2d}{d+1}\in (1, 2)
$$
and so $Q^+$ has a regularity-gaining property for
$\,\lambda\in \left(\frac d2, \,\frac {d+1}{2}\right]\,$
if $\,\alpha> d/4\,.$

Our discussions may be summarized as follows where we put
\begin{equation}\label{NN4}
\alpha_+(\lambda) = \left\{\aligned &{\,\,\,\,\frac d4\,} &{\quad\text{for}\quad \lambda\in \left(\frac 12\,,\,\frac{d+1}{2}\right]\,, }\\
& {\lambda - \frac d2} &{\quad\text{for}\quad \lambda\in\left(\frac {d+1}{2}\,,\,d\right)\,.}
\endaligned\right.
\end{equation}

\medskip

\begin{corollary}\label{corollaryN1} Assume that $\,B = |v-v_*|^{-\lambda}\,b\left(\mathbf{k}\cdot\sigma\right)\,$
where $\,b\in L^p(\mathbb{S}^{d-1})\,$ for some $\,2\le p\le\infty\,$ and  $\,1/2<\lambda<d\,.$
Let $\,f, g\in X^\alpha(\R^d)\,$ with $\,\alpha>\alpha_+(\lambda)\,.$
\begin{itemize}
\item[{\rm(i)}] There exists a constant $C$ independent of $f, g$ such that
\begin{equation}\label{NN5}
\left\|Q^+(f,g)\right\|_{\dot{H}^{\alpha + \delta}}\,\le\,C\,\|b\|_{L^p(\mathbb{S}^{d-1})}\,
\|f\|_{X^\alpha}\|g\|_{X^\alpha}
\end{equation}
for some $\,\delta = \delta(\lambda, p, d)\in [0, \,d-\lambda].\,$
\item[{\rm(ii)}] There exists a constant $C_+$ independent of $f, g$ such that
\begin{align}
\left\|Q^+(f,g)\right\|_{\dot{H}^{\alpha}}\,\le\,C_+\,\|b\|_{L^p(\mathbb{S}^{d-1})}\,
\|f\|_{X^\alpha}\|g\|_{X^\alpha}\,.\label{NN9}
\end{align}
\end{itemize}
\end{corollary}

\smallskip
\begin{proof} (i) is a summary of what we have shown in the above and (ii) follows instantly
upon combining Theorem \ref{theoremI1}
and (i) with the aid of Lemma \ref{lemmaHY3}.
\end{proof}

\medskip

\begin{remark}\label{remarkR2} In the case of $\,b\in L^p(\s^{d-1})\,$ with
$\,1<p<2,\,$ if we follow the foregoing analysis carefully, then it is not hard to
see that the maximum possible range of $\lambda$ is given by
\begin{equation*}\label{NN6}
\bigcup_{1<p_1\le p}\,I(p_1) = \left(\Phi(p),\,d\right)\,
\end{equation*}
and there exists $\,1<p_1<p\,$ satisfying $\,\left(\frac dp, \,1+ \frac {d-1}{p}\right]\subset I(p_1).\,$
Inspecting the minimally required condition on $\alpha$, hence, we conclude that
$Q^+$ has a regularity-gaining or preserving property on $X^\alpha(\R^d)$ in the sense of
(\ref{NN5}) for
$\,\lambda\in \left(\Phi(p),\,d\right)\,$ as long as
\begin{align}\label{NN7}
\alpha>\max\left\{ \frac d4, \, 1-\frac 1p +\left(\frac 1p -\frac 12\right)d,\,\lambda-\frac d2\right\}\,.
\end{align}
Of course, the corresponding estimate of type (\ref{NN9}) is also available.
\end{remark}

\subsection{Regularity of $Q^-$}
We now consider the regularity of the loss part $Q^-$ defined as
$$ Q^-(f,g)(v) = \|b\|_{L^1(\s^{d-1})}\,f(v)\,\int_{\R^d} |v- v_*|^{-\lambda} g(v_*)\,dv_*\,.
$$
If we make use of the Fourier representation formula (\ref{F3}) for the fractional integral of $g$,
it is straightforward to compute its Fourier transform
\begin{equation}\label{R0}
[Q^-(f,g)]\,\widehat\,\,(\xi) = \lambda_d\,\|b\|_{L^1(\s^{d-1})}\int_{\R^d} |\zeta|^{-d +\lambda}\,\hat{f}(\xi -\zeta)\,\hat{g}(\zeta)\,d\zeta\,.
\end{equation}

With the aid of Hausdorff-Young type inequalities, we prove that $Q^-$ has regularity-preserving properties in the following sense.

\medskip

\begin{theorem}\label{theoremR3} Assume that $\,B = |v-v_*|^{-\lambda}\,b\left(\mathbf{k}\cdot\sigma\right)\,$
where $\,b\in L^1(\mathbb{S}^{d-1})\,$ and $\,0<\lambda<d\,$. Let $\,f, g\,\in\,L^1\cap\dot{H}^\alpha(\R^d).\,$
\begin{itemize}
\item[{\rm(i)}] If $\,0<\lambda<d/2\,$ and $\,\alpha\ge 0,\,$ then
\begin{align}\label{R1}
&\left\|Q^-(f,g)\right\|_{\dot{H}^{\alpha}} \,\le\, C\,\|b\|_{L^1(\mathbb{S}^{d-1})}\nonumber\\
&\qquad\times\quad\biggl\{\|f\|_{\dot{H}^\alpha} \|g\|_{L^1}^{1-\theta}\|g\|_{\dot{H}^\alpha}^\theta
+\|f\|_{L^1}^{1-\delta}\|f\|_{\dot{H}^\alpha}^\delta\|g\|_{L^1}^{1-\gamma}\|g\|_{\dot{H}^\alpha}^{\gamma} \biggr\}
\end{align}
where $C$ is a constant depending only on $\,\alpha, \lambda, d\,$ and
$$\theta = \frac{2\lambda}{d+2\alpha},\quad\delta = \frac{d}{d+2\alpha}, \quad \gamma = \frac{2\lambda + 2\alpha}{d+ 2\alpha}.$$
\item[{\rm(ii)}] If $\,0<\lambda<d\,$ and $\,\alpha>d/2,\,$ then
there exists a constant $C$ depending only on $\,\alpha, \lambda, d\,$ such that
\begin{align}\label{R2}
&\left\|Q^-(f,g)\right\|_{\dot{H}^{\alpha}}\,\le\, C\,\|b\|_{L^1(\mathbb{S}^{d-1})}\nonumber\\
&\qquad\times\quad\biggl\{\|f\|_{\dot{H}^\alpha} \|g\|_{L^1}^{1-\theta}\|g\|_{\dot{H}^\alpha}^\theta
+\|f\|_{L^1}^{1-\mu}\|f\|_{\dot{H}^\alpha}^\mu\|g\|_{L^1}^{1-\nu}\|g\|_{\dot{H}^\alpha}^{\nu} \biggr\}
\end{align}
where $C$ is a constant depending only on $\,\alpha, \lambda, d\,$ and
$\,\mu = 2\delta, \,\, \nu = \gamma - \delta\,$ with the same $\,\theta, \delta, \gamma\,$ as above.
\end{itemize}
\end{theorem}

\smallskip

\begin{proof}
By using the elementary inequality $\,|\xi|^\alpha\le 2^\alpha\left(|\xi -\zeta|^\alpha + |\zeta|^\alpha\right)\,$
valid for any $\,\zeta\in\R^d\,$ and setting $\,\widehat{f_\alpha}(\xi) = |\xi|^\alpha \hat f(\xi)\,,\,\widehat{g_\alpha}(\xi) = |\xi|^\alpha \hat g(\xi),\,$ we estimate
\begin{align*}
\left\|Q^-(f,g)\right\|_{\dot{H}^{\alpha}}\,&\le\, \lambda_d\,2^\alpha\,\|b\|_{L^1(\s^{d-1})}\,\bigl\{{\rm(I)} + {\rm (II)}\bigr\}\quad\text{where}\\
 {\rm(I)}&= \left\|\int_{\R^d} |\zeta|^{-d +\lambda}\,\bigl|\widehat{f_\alpha}(\xi -\zeta)\bigr|\,\bigl|\hat{g}(\zeta)\bigr|\,d\zeta\right\|_{L^2(d\xi)}\,,\\
{\rm(II)}&= \left\|\int_{\R^d} |\zeta|^{-d +\lambda}\,\bigl|\hat{f}(\xi -\zeta)\bigr|\,\bigl|\widehat{g_\alpha}(\zeta)\bigr|\,d\zeta\right\|_{L^2(d\xi)}\,.
\end{align*}

For the first term, it follows from the Minkowski integral inequality and the estimate (\ref{F-3}) that if $\alpha$
satisfies the condition that $\,\alpha\ge 0\,$ when $\,0<\lambda<d/2\,$ and $\,\alpha>\lambda - d/2\,$ when $\,d/2\le\lambda<d,\,$ then
\begin{align}\label{R5}
{\rm(I)}&\le \bigl\|\widehat{f_\alpha}\bigr\|_{L^2}\,\int_{\R^d} |\zeta|^{-d +\lambda}\,\bigl|\hat{g}(\zeta)\bigr|\,d\zeta\nonumber\\
&\le C\, \|f\|_{\dot{H}^\alpha} \|g\|_{L^1}^{1-\theta}\|g\|_{\dot{H}^\alpha}^\theta\,,
\end{align}
where $C$ is a constant depending only on $\,\alpha, \lambda, d\,.$

For the second term, let us begin with the case $\,0<\lambda<d/2\,$ and $\,\alpha\ge 0\,.$ By applying the Minkowski integral inequality once again, we estimate
$${\rm(II)}\le \bigl\|\hat{f}\bigr\|_{L^2}\,\int_{\R^d} |\zeta|^{-d +\lambda}\,\bigl|\widehat{g_\alpha}(\zeta)\bigr|\,d\zeta\,.$$
For any $\,\rho>0\,,$ we observe
\begin{align*}
\int_{\R^d} |\zeta|^{-d +\lambda}\,\bigl|\widehat{g_\alpha}(\zeta)\bigr|\,d\zeta &\le C\,\left\{\rho^{\lambda +\alpha}\, \|g\|_{L^1} + \rho^{\lambda-d/2}\,\|g\|_{\dot{H}^\alpha}\right\}\,,
\end{align*}
a simple consequence of splitting the integral into
two parts by decomposing $\,\R^d = \{|\zeta|\le\rho\}\cup \{|\zeta| >\rho\}\,$
and applying the Cauchy-Schwartz inequality for the latter part. Optimizing in $\rho$ gives
$$\int_{\R^d} |\zeta|^{-d +\lambda}\,\bigl|\widehat{g_\alpha}(\zeta)\bigr|\,d\zeta\,\le C\, \|g\|_{L^1}^{1-\gamma}\|g\|_{\dot{H}^\alpha}^{\gamma}$$
with the same $\gamma$ defined as in the statement of (i). Combining with the Hausdorff-Young type inequality (\ref{HY2}) for $\,\|\hat f\|_{L^2}\,,$
we get
\begin{equation}\label{R6}
{\rm(II)}\,\le \, C\, \biggl\{\|f\|_{L^1}^{1-\delta}\|f\|_{\dot{H}^\alpha}^\delta\|g\|_{L^1}^{1-\gamma}\|g\|_{\dot{H}^\alpha}^{\gamma} \biggr\}\,.
\end{equation}
Adding the estimates (\ref{R5}) and (\ref{R6}), we complete the proof of part (i).

In the case $\,0<\lambda<d\,$ and $\,\alpha>d/2\,,$ we change variables repeatedly and apply the Minkowski integral inequality
to observe
\begin{align*}
{\rm(II)}&= \left\|\int_{\R^d} |\xi-\zeta|^{-d +\lambda}\,\bigl|\widehat{g_\alpha}(\xi-\zeta)\bigr|\,\bigl|\hat{f}(\zeta)\bigr|\,d\zeta\right\|_{L^2(d\xi)}\\
&\le \bigl\|\hat{f}\bigr\|_{L^1}\,\left(\int_{\R^d} |\xi|^{-2d +2\lambda}\,\bigl|\widehat{g_\alpha}(\xi)\bigr|^2\,d\xi\right)^{1/2}\,.
\end{align*}
Since $\,-d + 2\lambda + 2\alpha>0\,,$ it is easy to see
$$ \int_{\R^d} |\xi|^{-2d +2\lambda}\,\bigl|\widehat{g_\alpha}(\xi)\bigr|^2\,d\xi\,\le\,C\,
\left\{\rho^{-d + 2\lambda +2\alpha}\, \|g\|_{L^1}^2 + \rho^{-2d + 2\lambda}\,\|g\|_{\dot{H}^\alpha}^2\right\}
$$
for any $\,\rho>0\,$ and optimizing in $\rho$ gives
$$ \int_{\R^d} |\xi|^{-2d +2\lambda}\,\bigl|\widehat{g_\alpha}(\xi)\bigr|^2\,d\xi\,\le\,C\,
\|g\|_{L^1}^{2(1-\nu)}\|g\|_{\dot{H}^\alpha}^{2\nu}\,.$$
Combining with the Hausdorff-Young type inequality (\ref{HY2}) for $\,\|\hat f\|_{L^1}\,$ valid
due to the condition $\,\alpha>d/2\,,$ we finally obtain
\begin{equation}\label{R7}
{\rm(II)}\,\le \, C\, \biggl\{\|f\|_{L^1}^{1-\mu}\|f\|_{\dot{H}^\alpha}^\mu\|g\|_{L^1}^{1-\nu}\|g\|_{\dot{H}^\alpha}^{\nu} \biggr\}\,.
\end{equation}
Adding the estimates (\ref{R5}) and (\ref{R7}), we complete the proof of part (ii).
\end{proof}

\medskip

An immediate consequence is the following:
\medskip

\begin{corollary}\label{corollaryR1}
Assume that $\,B = |v-v_*|^{-\lambda}\,b\left(\mathbf{k}\cdot\sigma\right)\,$
where $\,b\in L^1(\mathbb{S}^{d-1})\,$ and $\,0<\lambda<d\,$.
If $\,f, g\in X^\alpha(\R^d)\,$ with $\,\alpha\ge 0\,$ for $\,0<\lambda<d/2\,$ and $\,\alpha> d/2\,$ for $\,d/2\le\lambda<d,\,$
there exists a constant $C_{-}$ independent of $f, g$ such that
\begin{align}
\left\|Q^-(f,g)\right\|_{\dot{H}^{\alpha}}\,\le\,C_-\,\|b\|_{L^1(\mathbb{S}^{d-1})}\,
\|f\|_{X^\alpha}\|g\|_{X^\alpha}\,.\label{R4}
\end{align}
\end{corollary}

\medskip
\subsection{Regularity Summary}
If we combine Corollary \ref{corollaryN1} and Corollary \ref{corollaryR1},
we obtain what we aim to establish for the regularity of $Q$.

\medskip

\begin{theorem}\label{theoremQR1}
Assume that $\,B = |v-v_*|^{-\lambda}\,b\left(\mathbf{k}\cdot\sigma\right)\,$ where
$\,b\in L^p(\mathbb{S}^{d-1})\,$ for some $\,2\le p\le\infty\,$ and  $\,1/2<\lambda<d\,.$
If $\,f, g\in X^\alpha(\R^d)\,$ with $\alpha$ satisfying
$\,\alpha>d/4\,$ for $\,1/2<\lambda<d/2\,$ and $\,\alpha> d/2\,$ for $\,d/2\le\lambda<d,\,$
then there exists a constant $C_{R}$ independent of $f, g$ such that
\begin{align}
\left\|Q(f,g)\right\|_{\dot{H}^{\alpha}}\,\le\,C_R\,\|b\|_{L^p(\mathbb{S}^{d-1})}\,
\|f\|_{X^\alpha}\|g\|_{X^\alpha}\,.\label{QR1}
\end{align}
Moreover, the quadratic operator $\,f\mapsto Q(f,f)\,$ satisfies
\begin{align}
&\left\|Q(f,f)\right\|_{\dot{H}^{\alpha}}\,\le\,C_R\,\|b\|_{L^p(\mathbb{S}^{d-1})}\,
\|f\|_{L^1}^{1-\theta}\|f\|_{\dot{H}^\alpha}^{1+\theta}\label{QR2}\\
&\qquad\text{with}\quad \theta = \frac{2\lambda}{d+2\alpha}\,\in\,(0, 1).\nonumber
\end{align}
\end{theorem}

\smallskip

\begin{proof}
Keeping track of multiplicative constants, if we put
$$\,C_R = C_+ + C_{-}\,|\s^{d-1}|^{1-1/p}\,$$ where $\,C+, \,C_{-}\,$ are the same positive constants
appearing in the $Q^+$-regularity estimate (\ref{NN9}) and the $Q^-$-regularity estimates (\ref{R4}), respectively,
then the stated bilinear estimate (\ref{QR1}) follows at once.

As for the quadratic estimate (\ref{QR2}), if we inspect Theorem \ref{theoremR2} and the reasonings leading to
Corollary \ref{corollaryN1} carefully, it is not hard to find
\begin{align}\label{NNN5}
&\left\|Q^+(f,f)\right\|_{\dot{H}^{\alpha + \delta}}\,\le\,C\,\|b\|_{L^p(\mathbb{S}^{d-1})}\,
\|f\|_{L^1}^{1-\nu}\|f\|_{\dot{H}^\alpha}^{1+\nu}\\
&\text{with}\quad \nu = \frac{2(\delta +\lambda)}{d+2\alpha}\,\in\,(0, 1)\,,\nonumber
\end{align}
under the same setting as in Corollary \ref{corollaryN1}, which in turn implies
\begin{align}
\left\|Q^+(f,f)\right\|_{\dot{H}^{\alpha}}\,\le\,C_+\,\|b\|_{L^p(\mathbb{S}^{d-1})}\,
\|f\|_{L^1}^{1-\theta}\|f\|_{\dot{H}^\alpha}^{1+\theta}\label{NNN9}\,.
\end{align}
On the other hand, it is straightforward to find
\begin{align}
\left\|Q^-(f,f)\right\|_{\dot{H}^{\alpha}}\,\le\,C_-\,\|b\|_{L^1(\mathbb{S}^{d-1})}\,
\|f\|_{L^1}^{1-\theta}\|f\|_{\dot{H}^\alpha}^{1+\theta}\label{RR4}
\end{align}
under the same setting as in Corollary \ref{corollaryR1}. Adding these two estimates (\ref{NNN9}) and (\ref{RR4}),
we obtain the desired estimate (\ref{QR2}).

\end{proof}

\medskip
\begin{remark} While these results will be our final forms of regularity to be used in the subsequent sections,
we point out that the range of $\,\lambda\in (0, 1/2]\,$ is missing due to lack of
regularity estimates for the gain part $Q^+$ as in Corollary \ref{corollaryN1}.
If $f, g$ are assumed to satisfy extra condition $\,vf, vg\in \dot{H}^\alpha(\R^d)\,$ with $\,\alpha>\lambda\,,$
however, it can be shown that $Q^+$ gains regularity of order $(d-1)/2$
(for the details and ideas, see Theorem 5.1 of \cite{Cho} and \cite{BD}).

In comparison with our work, we also refer to the paper \cite{MV} of C. Mouhot and
C. Villani  for the regularity theory in the case of cutoff hard potentials for which
a notable difference is that certain lower bounds are available for the regularity estimates
of the loss part $Q^-$.
\end{remark}

\section{Proofs of the Main Results}
Having obtained bilinear estimates for the integrability and regularity of the collision
operator $Q$ with cutoff soft potentials, we are now in a position to prove our main results.

\paragraph{Proof of Theorem 1.}
Under the stated hypotheses on $\,b, \lambda\,$ and
$\alpha$, we may rephrase the integrability estimate (\ref{I2}) of Corollary \ref{corollaryI1} and
the regularity estimate (\ref{QR1}) of Theorem \ref{theoremQR1} in the unifying form
\begin{equation}\label{P1}
\left\| Q(f, g)\right\|_{X^\alpha}\,\le\,\max\,\bigl(C_I,\,C_L\bigr)\,\|b\|_{L^p(\mathbb{S}^{d-1})}\,\|f\|_{X^\alpha}\,\|g\|_{X^\alpha}
\end{equation}
for all $\,f, g\in X^\alpha(\R^d)\,.$ We define the time $T$ in question by
\begin{equation}\label{P2}
T = \frac{1}{\,\,\,5 \max\bigl(C_I,\,C_R\bigr)\,\|b\|_{L^p(\mathbb{S}^{d-1})}\,\|f_0\|_{X^\alpha}\,\,}
\end{equation}
where both $\,\|b\|_{L^p(\mathbb{S}^{d-1})}\,$ and $\,\|f_0\|_{X^\alpha}\,$ are assumed to be positive
without any loss of generality. With this $T$, we proceed to prove Theorem 1 along
three categories. For simplicity, let us put
\begin{equation}\label{PK}
K_b = \max\,\bigl(C_I,\,C_L\bigr)\,\|b\|_{L^p(\mathbb{S}^{d-1})}\,.
\end{equation}

\subsection{Existence and Uniqueness} Let $\mathcal{A}$ denote the operator defined by
\begin{equation*}
\left(\mathcal{A} f\right)(v, t) = f_0(v) + \int_0^t Q(f, f)(v, s) \,ds\quad\text{for}\quad t\in [0, T].
\end{equation*}
We prove the existence and uniqueness part of our theorem by showing that $\mathcal{A}$ has a unique fixed point in the
space $\Omega_T$ defined as
\begin{equation*}
\Omega_T = \left\{ f\in C\left([0, T] ; X^\alpha(\R^d)\right)\,:\,
\sup_{t\in [0,\, T]}\,\left\| f(\cdot, t)\right\|_{X^\alpha}\le \,2\left\| f_0\right\|_{X^\alpha}\,\right\}\,,
\end{equation*}
which is a complete metric space with respect to the induced metric
\begin{equation}
d_T(f, g) = \sup_{t\in [0, \,T]}\,\left\| (f -g)(\cdot, t)\right\|_{X^\alpha}\,.
\end{equation}

\begin{itemize}
\item[(1)]
We first show that $\mathcal{A}$ maps $\Omega_T$ into itself.
In view of Fubini's theorem and Minkowski's integral inequality, it is trivial to see
\begin{align*}
\left\|\int_0^t Q(f, f)(\cdot, s) ds\right\|_{L^1} &\le \int_0^t \left\| Q(f, f)(\cdot, s)\right\|_{L^1}\,ds\,,\\
\left\|\int_0^t Q(f, f)(\cdot, s) ds\right\|_{\dot{H}^\alpha} &\le \int_0^t\left\| Q(f, f)(\cdot, s)\right\|_{\dot{H}^\alpha}\,ds\,.
\end{align*}
For $\,f\in\Omega_T\,$ and $\,t\in[0, T]\,,$ it follows from (\ref{P1}) that
\begin{align}\label{P3}
\left\| \left(\mathcal{A}f\right)(\cdot, t)\right\|_{X^\alpha} &\le \|f_0\|_{X^\alpha} +  \int_0^t\,\left\| Q(f, f)(\cdot, s)\right\|_{X^\alpha}\,ds\nonumber\\
&\le \|f_0\|_{X^\alpha} + K_b\,\int_0^t\,\|f(\cdot, s)\|_{X^\alpha}^2\,ds
\end{align}
whence it is trivial to observe
\begin{equation}\label{P5}
\sup_{t\in [0, \,T]}\,\left\| \left(\mathcal{A}f\right)(\cdot, t)\right\|_{X^\alpha}\,\le\, \frac 95\,\left\|f_0\right\|_{X^\alpha}\,.
\end{equation}

The $t$-continuity of $\,\|(Af)(\cdot, t)\|_{X^\alpha}\,$ follows from the evident estimate
\begin{equation}\label{P6}
\left\| \left(\mathcal{A}f\right)(\cdot, t) - \left(\mathcal{A}f\right)(\cdot, s)\right\|_{X^\alpha} \,
\le\, 4 K_b \left\|f_0\right\|_{X^\alpha}^2\,\left|\,t - s\,\right|
\end{equation}
for all $\, t, s\in [0, T].\,$ By (\ref{P5}), (\ref{P6}), $\,\mathcal{A} f\in\Omega_T\,.$

\item[(2)] We next show that $\mathcal{A}$ is a contraction mapping on $\Omega_T$.
Fix $\,f, g\in\Omega_T\,.$ Due to the bilinearity
of $Q$, we note that
$$Q(f, f) - Q(g, g) = \frac 12 \left[ Q(f-g, f+g) + Q(f+g, f-g)\right]\,.$$
For $\,t\in [0, T]\,,$ this identity and the estimate (\ref{P1}) imply
\begin{align*}
\left\| \left(\mathcal{A}f - \mathcal{A}g\right)(\cdot, t)\right\|_{X^\alpha} &= \left\|\int_0^t  \bigl[Q(f, f)(\cdot, s)-Q(g, g)(\cdot, s)\bigr]\,ds
\right\|_{X^\alpha}\\
&\le K_b\,\int_0^t\,\|(f-g)(\cdot, s)\|_{X^\alpha}\,\|(f+g)(\cdot, s)\|_{X^\alpha}\,ds\\
&\le \frac 45\,\left(\sup_{t\in[0, \,T]} \|(f-g)(\cdot, t)\|_{X^\alpha}\right)
\end{align*}
and thus taking supremum yields
\begin{equation}\label{P9}
d_T\left(\mathcal{A}f,\,\mathcal{A}g\right)\,\le\,\frac 45\,d_T(f, g)\,,
\end{equation}
which shows that $A$ is a contraction on $\Omega_T$.
\end{itemize}

By the contraction mapping principle, we conclude that $\mathcal{A}$ has a unique fixed point
$f$ in $\Omega_T$, that is, a unique function
$\,f\in C\left([0, T]; X^\alpha(\R^d)\right)\,$ solving
\begin{equation}\label{P12}
f(v, t) = f_0(v) + \int_0^t Q(f, f)(v, s)\,ds
\end{equation}
on $\,\R^d\times [0, T],\,$ as we wished to prove.

\subsection{Strong Differentiability}
That $f$ conserves mass on $\,[0, \,T]\,$ is obvious due to
$$\int_{\R^d} Q^+(f, f)(v, t)\,dv = \int_{\R^d} Q^-(f, f)(v, t)\,dv\,.$$

To prove part (ii), fix $\,t\in (0, T).\,$ For $\,\tau>0\,$ sufficiently small, (\ref{P12}) gives
\begin{align*}
&\frac{f(v, t+\tau) - f(v, t)}{\tau}  - Q(f, f)(v, t)\nonumber\\
&\qquad\qquad= \frac 1\tau\int_t^{t+\tau}
\left[Q(f, f)(v, s) - Q(f, f)(v, t)\right] ds\nonumber
\end{align*}
for any $\,v\in\R^d.\,$ Taking $X^\alpha$ norm, we deduce from (\ref{P1})
\begin{align}\label{P10}
&\left\|\frac{f(\cdot, t+\tau) - f(\cdot, t)}{\tau} - Q(f, f)(\cdot, t)\right\|_{X^\alpha}\nonumber\\
&\qquad\le \frac 1\tau\int_t^{t+\tau}\left\|Q(f, f)(\cdot, s) - Q(f, f)(\cdot, t)\right\|_{X^\alpha}\, ds\nonumber\\
&\qquad\le \frac {\,K_b\,}{\,\,\tau\,\,}\int_t^{t+\tau} \|f(\cdot, s) - f(\cdot, t)\|_{X^\alpha}\|f(\cdot, s) + f(\cdot, t)\|_{X^\alpha}\,ds\nonumber\\
&\qquad\le \left(16 K_b^2 \|f_0\|_{X^\alpha}^3\right)\cdot\tau\nonumber
\end{align}
where we used the continuity estimate (\ref{P6}) and the fact $\,f\in\Omega_T\,.$
By carrying out the same analysis for $\,\tau<0,\,$ we get
\begin{equation}\label{P11}
\left\|\frac{f(\cdot, t+\tau) - f(\cdot, t)}{\tau} - Q(f, f)(\cdot, t)\right\|_{X^\alpha}
\le \left(16 K_b^2 \|f_0\|_{X^\alpha}^3\right)\cdot |\tau|
\end{equation}
for all $\,\tau\ne 0\,$ sufficiently small. Letting $\,\tau\to 0,\,$ we conclude
$$\frac{\partial f}{\partial t} = Q(f, f) \quad\text{in}\quad X^\alpha(\R^d)\quad\text{for}\quad t\in (0, T).$$

The same arguments give
\begin{equation}\label{P13}
\left\|Q(f, f)(\cdot, t) - Q(f, f)(\cdot, s)\right\|_{X^\alpha}\,\le\,\left(16 K_b^2 \|f_0\|_{X^\alpha}^3\right)\cdot |t-s|
\end{equation}
for all $\,t, s\in [0, T]\,,$ which proves the t-continuity of $\,\left\|\partial_t f(\cdot, t)\right\|_{X^\alpha}\,$
and the proof of part (ii) is complete.

\subsection{Nonnegativity.}
It remains to prove that the unique solution $f$ constructed as above is nonnegative. To do this, we shall modify the
idea of X. Lu and Y. Zhang (\cite{LZ}) which originally deals with the nonnegativity problem
for solutions to the inhomogeneous Boltzmann
equation in $\R^3$. Let us put
\begin{equation}\label{P14}
N_f(t) = \int_{\R^d}\left( -f(v, t)\right)^+\,dv\qquad(t\ge 0)\,,
\end{equation}
where $\,(-f)^+ = \max\left( -f, 0\right), \,$ that is, the negative part of $f$. Let $\chi$ denote the characteristic
function of the interval $(0, \infty)$. Since $f_0$ is nonnegative,
\begin{align}\label{P15}
N_f(t) &= \int_0^t\int_{\R^d}\bigl[-Q(f,f)(v, s)\bigr]\,\chi\left(-f(v, s)\right)\,ds\nonumber\\
&=\int_0^t\int_{\R^d}\int_{\R^d}\int_{\s^{d-1}} B\,\left( -f' f_*' + f f_*\right)\,\chi(-f)\,d\sigma dv_* dv ds
\end{align}
with the usual abbreviation notation for $\,f', f_*'\,.$

If we use the inequality
\begin{equation}\label{P16}
\left( -f' f_*' + f f_*\right)\,\chi(-f)\,\le\, \left(-f' f_*'\right)^+ - \left(-f f_*\right)^+ + |f|(- f_*)^+\,,
\end{equation}
which is proved in Lemma 2 of \cite{LZ}, then we deduce from (\ref{P15})
\begin{align}\label{P17}
N_f(t) & \le \int_0^t\int_{\R^d}\int_{\R^d}\int_{\s^{d-1}} B\,\left[\left(-f' f_*'\right)^+ - \left(-f f_*\right)^+\right]
\,d\sigma dv_* dv ds\nonumber\\
& \qquad + \, \int_0^t\int_{\R^d}\int_{\R^d}\int_{\s^{d-1}} B\,|f|(- f_*)^+\, d\sigma dv_* dv ds\nonumber\\
&=\int_0^t\int_{\R^d}\int_{\R^d}\int_{\s^{d-1}} B\,|f|(- f_*)^+\, d\sigma dv_* dv ds
\end{align}
where the first term vanishes due to
$$ \int_{\R^d}\int_{\R^d}\int_{\s^{d-1}} B\,\left(-f' f_*'\right)^+ d\sigma dv_* dv
= \int_{\R^d}\int_{\R^d}\int_{\s^{d-1}} B\,\left(-f f_*\right)^+ d\sigma dv_* dv\,.$$

By making use of the estimate (\ref{F2}) for the functional  $J_\lambda(f)$ established in Lemma \ref{lemmaF1},
we deduce from (\ref{P17})
\begin{align}\label{P18}
N_f(t) &\le \|b\|_{L^1(\s^{d-1})} \,\int_0^t\int_{\R^d}\int_{\R^d} |v-v_*|^{-\lambda} |f(v, s)|(- f(v_*, s))^+\,dv_* dv ds \nonumber\\
&\le \|b\|_{L^1(\s^{d-1})} \,\int_0^t\, J_\lambda(f)(s)\, \left[\int_{\R^d}(- f(v_*, s))^+\,dv_*\right]\,ds\nonumber\\
&\le C_J\,\|b\|_{L^1(\s^{d-1})} \,\int_0^t\,\| f(\cdot, s)\|_{X^\alpha}\,N_f(s)\, ds\,.
\end{align}
If we recall $\,C_I = 2C_J\,|\s^{d-1}|^{1-1/p}\,,$ then (\ref{P18}) gives
\begin{align*}
N_f(t)\le C_I\,\|b\|_{L^p(\s^{d-1})}\,\|f_0\|_{X^\alpha}\,\int_0^t N_f(s)\,ds \le \frac{1}{5}\, \left(\sup_{t\in [0,\,T]}\,N_f(t)\right)
\end{align*}
for all $\,t\in [0, T]\,.$ Consequently,
\begin{equation}
\sup_{t\in [0,\,T]}\,N_f(t)\,\le\, \frac{1}{5}\, \left(\sup_{t\in [0,\,T]}\,N_f(t)\right)\,,
\end{equation}
which implies at once $\, N_f(t) = 0 \,$ for all $\,t\in [0, T]\,.$

An immediate conclusion is that $\,f\ge 0\,$ a.e. on $\,\R^d\times [0, T]\,.$
Using $\,f= |f|\,$ a.e. and arguing as in X. Lu and Y. Zhang (\cite{LZ}), it can be shown that $|f|$ is also a solution to
the integral equation (\ref{P12}) in $\Omega_T$. By the uniqueness, we conclude that $f$ is nonnegative everywhere
on the interval $\,[0, \,T]\,.$

Our proof of Theorem 1 is now complete.

\medskip

\paragraph{Proof of Corollary 1.} By the Sobolev embedding $\,X^\alpha(\R^d)\subset C_0(\R^d)\,$ for $\,\alpha>d/2\,$
as in Lemma \ref{lemmaSE}, the estimate (\ref{P11}) implies
\begin{equation}
\sup_{v\in\R^d}\,\left|\frac{f(v, t+\tau) - f(v, t)}{\tau} - Q(f, f)(v, t)\right|\le C|\tau|
\end{equation}
for sufficiently small $\,\tau\ne 0\,$ and for some constant $C$ independent of $\,\tau, t\,,$
which is enough to conclude $f$ is a classical solution to the Boltzmann equation (\ref{1.1}). The other stated
properties are also trivial consequences of this Sobolev embedding.

\section{Concluding Remarks}
In the above construction, the time interval $\,[0, T]\,$ for which the unique solution $f$ exists
on $\,\R^d\times [0, T]\,$ is determined by
$$ T= \frac{1}{\,5 K_b\,\|f_0\|_{X^\alpha}\,}$$
where $K_b$ denotes the multiplicative constant of the bilinear estimate (\ref{P1}) defined as in (\ref{PK}).
Since $f$ satisfies
\begin{align*}
\|f(\cdot, t)\|_{X^\alpha} &\le \|f_0\|_{X^\alpha} + \int_0^t\left\| Q(f, f)(\cdot, s)\right\|_{X^\alpha}\,ds\\
&\le \|f_0\|_{X^\alpha} + K_b\,\int_0^t\left\|f(\cdot, s)\right\|_{X^\alpha}^2\,ds\,,
\end{align*}
an easy modification of Gronwall's lemma yields
\begin{align}\label{CR2}
\|f(\cdot, t)\|_{X^\alpha} \le \frac{\|f_0\|_{X^\alpha}}{\,1- K_b\,\|f_0\|_{X^\alpha}\,t\,}
\le \frac 54\,\|f_0\|_{X^\alpha}\quad\text{for all}\quad t\in [0, T]\,.
\end{align}

In an attempt to extend the time interval of existence, let us put
\begin{equation}\label{CR3}
T_1 = T + \frac{1}{\,5K_b\,\|f(\cdot, T)\|_{X^\alpha}\,}\,.
\end{equation}
On the time interval $\,[T, T_1],\,$ if we choose $f(v, T)$ as the new initial datum and
repeat the same arguments of the contraction mapping principle as well as nonnegativity, then we obtain a unique solution $f_1$ to
\begin{equation*}
f_1(v, t) = f(v, T) + \int_T^t Q(f_1, f_1)(v, s)\,ds
\end{equation*}
on $\,\R^d\times [T, T_1]\,$ which satisfies
\begin{align}
\left\|f_1(\cdot, t)\right\|_{X^\alpha} \le \frac{\|f(\cdot, T)\|_{X^\alpha}}{\,1- K_b\|f(\cdot, T)\|_{X^\alpha}\,(t-T)\,}
\le \frac 54\,\|f(\cdot, T)\|_{X^\alpha}
\end{align}
for all $\, t\in [T, T_1].\,$ Gluing the two solutions together, if we define
\begin{equation*}
\tilde{f}(v, t) = \left\{\aligned &{f(v, t)\,\,\quad\text{for}\quad (v, t)\in \R^d\times [0, T], }\\
&f_1(v, t)\quad\text{for}\quad (v, t)\in \R^d\times [T, T_1],
\endaligned\right.
\end{equation*}
then it is obvious that $\tilde{f}$ is a unique solution to
\begin{equation*}
\tilde{f}(v, t) = f_0(v) + \int_0^t Q(\tilde{f}, \tilde{f})(v, s)\,ds
\end{equation*}
on $\,\R^d\times [0, T_1]\,.$ For simplicity, let us denote this extended solution $\tilde{f}$ by $f$.

Continuing this process, if $\,T<T_1<\cdots<T_\ell\,$ were selected and a unique solution $f$
were obtained on the time interval $\,[0, T_\ell]\,,$ we put
\begin{align}\label{CR4}
T_{\ell +1} = T_\ell + \frac{1}{\,5K_b\,\|f(\cdot, T_\ell)\|_{X^\alpha}\,}
\end{align}
and extend $f$ to $\,[0, T_{\ell +1}]\,$ in such a unique way that it solves
\begin{align*}
f(v, t) = f_0(v)  + \int_{0}^t Q(f, f)(v, s)\,ds
\end{align*}
on $\,\R^d\times [0, T_{\ell+1}]\,$ and satisfies the norm estimate
\begin{align}\label{CR5}
\|f(\cdot, t)\|_{X^\alpha} \le \frac 54 \|f(\cdot, T_\ell)\|_{X^\alpha}\quad\text{for all}\quad t\in [T_\ell, T_{\ell+1}]\,.
\end{align}

This iterative scheme, repeated infinitely many times if necessary, may produce a unique solution $f$
in $X^\alpha(\R^d)$ to the Boltzmann equation (\ref{1.1}) either for all time or for finite time.
While the global-in-time existence is inconclusive, (\ref{CR4}) and (\ref{CR5}) imply
\begin{equation*}
T_{\ell +1} - T_\ell \ge \frac 45\,\left( T_\ell - T_{\ell -1}\right)
\end{equation*}
for all $\,\ell\ge 1\,,$ where we put $\,T_0=T\,,$ and consequently
$$ T + \sum_{\ell\ge 1} \left(T_{\ell +1} - T_\ell\right) \,\ge\, T\,\sum_{\ell\ge 0}\left(\frac 45\right)^\ell =  5T = \frac{1}{\,K_b\,\|f_0\|_{X^\alpha}\,}\,.$$
Thus this iterative scheme covers at least the time interval $\,[0, \,5T)\,.$

It is now clear that there exists a maximal time of existence $T_*$, possibly
infinite, such that there exists a unique solution $f$ in $X^\alpha(\R^d)$ to the Boltzmann equation (\ref{1.1})
on $\,\R^d\times [0, T_*)\,$. If it is finite, then necessarily
\begin{equation}\label{CR6}
\lim_{t\,\nearrow\, T_*} \left\| f(\cdot, t)\right\|_{X^\alpha} = +\infty
\end{equation}
for otherwise we may take $\,f(v, T_*)\,$ as the new initial datum and repeat the same process to extend $f$
to a solution on a larger time interval. Due to the conservation of mass for all time,
the necessary condition (\ref{CR6}) may be considered in terms of the Sobolev norm $\,\|\cdot\|_{\dot{H}^\alpha}\,$.

As a summary, we state the following where we write
\begin{align}
\|Q\|_{X^\alpha\to X^\alpha} = \inf\,\biggl\{\,C>0\,:\, &\left\| Q(f, g)\right\|_{X^\alpha}
\le C\|b\|_{L^p(\s^{d-1})}\|f\|_{X^\alpha} \|g\|_{X^\alpha}\,\nonumber\\
&\qquad \text{for all}\quad f, g\in X^\alpha(\R^d)\,\biggr\}\,,
\end{align}
the bilinear operator norm of $Q$ on the pair of $X^\alpha(\R^d)$.

\medskip

\begin{corollary}
Under the same hypotheses as in Theorem 1, there exists a maximal time of existence $T_*$ with
\begin{equation}\label{CR7}
\frac{1}{\,\,\|Q\|_{X^\alpha\to X^\alpha}\,\|b\|_{L^p(\s^{d-1})}\,\|f_0\|_{X^\alpha}\,\,}\,\le\, T_*\,\le\,\infty\,
\end{equation}
for which the Boltzmann equation (\ref{1.1}) has a unique solution $f$ in $X^\alpha(\R^d)$
on $\,\R^d\times [0, T_*)\,$. Moreover, if $\,T_*<\infty\,,$ then necessarily
\begin{equation}\label{CR8}
\lim_{t\,\nearrow\, T_*} \left\| f(\cdot, t)\right\|_{\dot{H}^\alpha} = +\infty\,.
\end{equation}
\end{corollary}

\medskip

\begin{remark}
A slightly better lower bound of $T_*$ could be obtained if we make use of the
Sobolev regularity estimate (\ref{QR2}) for the quadratic term and apply a Bernoulli-type Gronwall's lemma
in view of the conservation of mass for all time. Indeed, it is easy to see
\begin{equation}\label{CR9}
\|f(\cdot, t)\|_{\dot{H}^\alpha}\le \frac{\|f_0\|_{\dot{H}^\alpha}}
{\,\,\left( 1- \theta\,C_R\|b\|_{L^p(\s^{d-1})}\,\|f_0\|_{L^1}^{1-\theta}\,\|f_0\|_{\dot{H}^\alpha}^\theta\,t\right)^{1/\theta}\,\,}
\end{equation}
for all $\,t\in[0, T]\,$ where $\,\theta = 2\lambda/(d+2\alpha)\in (0, 1)\,$ and the same iterative extension as above yields
a lower bound of $T_*$.
\end{remark}

Such a refinement, however, still does not yield any conclusive information on the global existence or on the finite
time condition (\ref{CR8}). Except a definite lower bound, the scenario of global existence of a smooth solution
to the Boltzmann equation (\ref{1.1}) with soft potentials is similar to those of many other evolutionary
partial differential equations such as the Euler or Navier-Stokes equations (see Corollary 3.2, \cite{BM}, for instance).

\medskip

As our final comment, let us point out that a different choice of $p$ in our $L^p$ cutoff conditions on $b$ influences the range of $\lambda$
and the required degree of smoothness $\alpha$ regarding the regularity properties of $Q^+$, which can be notably seen from
Corollary \ref{corollaryN1} and Remark \ref{remarkR2}. Although it is not clear what it means from a physical
point of view, the level of influence becomes more complicated as $p$ gets smaller being close to $1$. For simplicity, we only
state the following local existence theorem of a smooth classical solution.

\medskip

\begin{corollary} Assume that $\,B= |v-v_*|^{-\lambda} b(\mathbf{k}\cdot\sigma)\,$ where $\,b\in L^p(\s^{d-1})\,$ for some
$\,1<p <2\,$ and
\begin{equation}\label{CR10}
\frac dp - (d-1)\left(1-\frac 1p\right) <\lambda<d\,.
\end{equation}
Let $f_0$ be a nonnegative initial datum in $X^\alpha(\R^d)$ with $\,\alpha>d/2\,.$
Then there exists $\,T>0,\,$ depending only on $\,\alpha, \lambda, p, d, \,\|b\|_{L^p(\s^{d-1})}, \,\|f_0\|_{X^\alpha},\,$ such that
the Boltzmann equation (\ref{1.1}) has a unique classical solution $f$ with the same properties as stated in Theorem 1
and Corollary 1.
\end{corollary}

\end{rmfamily}

\end{document}